\documentclass[reqno]{amsart}

\usepackage{xypic}
\input xy
\xyoption{all}
\usepackage{epsfig}
\usepackage{color}
\usepackage{amsthm}
\usepackage{amssymb}
\usepackage{amsmath}
\usepackage{amscd}

\newcommand{\labell}[1] {\label{#1}}

\numberwithin{equation}{section}
\newtheorem {Theorem}{Theorem} 
\numberwithin{Theorem}{section}

\newtheorem {Proposition}[Theorem]{Proposition}  
\theoremstyle{definition}
\newtheorem{Definition}[Theorem]{Definition}
\theoremstyle{remark}
\newtheorem{Remark}[Theorem]{Remark}

\newtheorem {Corollary}[Theorem]{Corollary}  

\expandafter\chardef\csname pre amssym.def at\endcsname=\the\catcode`\@
\catcode`\@=11
\def\undefine#1{\let#1\undefined}
\def\newsymbol#1#2#3#4#5{\let\next@\relax
 \ifnum#2=\@ne\let\next@\msafam@\else
 \ifnum#2=\tw@\let\next@\msbfam@\fi\fi
 \mathchardef#1="#3\next@#4#5}
\def\mathhexbox@#1#2#3{\relax
 \ifmmode\mathpalette{}{\m@th\mathchar"#1#2#3}%
 \else\leavevmode\hbox{$\m@th\mathchar"#1#2#3$}\fi}
\def\hexnumber@#1{\ifcase#1 0\or 1\or 2\or 3\or 4\or 5\or 6\or 7\or 8\or
 9\or A\or B\or C\or D\or E\or F\fi}

\font\teneufm=eufm10
\font\seveneufm=eufm7
\font\fiveeufm=eufm5
\newfam\eufmfam
\textfont\eufmfam=\teneufm
\scriptfont\eufmfam=\seveneufm
\scriptscriptfont\eufmfam=\fiveeufm

\catcode`\@=\csname pre amssym.def at\endcsname

\def    \eps    {\epsilon}

\newcommand{\supp}{\operatorname{supp}}

\newcommand{\Ham}{{\mathit Ham}}

\newcommand{\Cc}{{\mathcal C}}

\newcommand{\Hh}{{\mathcal H}}

\newcommand{\Mm}{{\mathcal M}}

\newcommand{\Pp}{{\mathcal P}}

\newcommand{\Ss}{{\mathcal S}}

\def    \C      {{\mathbb C}}
\def    \R      {{\mathbb R}}

\def    \Z      {{\mathbb Z}}
\def    \N      {{\mathbb N}}

\def    \T      {{\mathbb T}}
\def    \CP     {{\mathbb C}{\mathbb P}}

\def    \12    {{\frac{1}{2}}}
\def    \p      {\partial}

\def    \HF     {\operatorname{HF}}
\def    \HM     {\operatorname{HM}}
\def    \HZ     {\operatorname{HZ}}

\def    \CF     {\operatorname{CF}}

\def    \Ham    {\operatorname{Ham}}

\def    \MUCZ  {\operatorname{\mu_{\scriptscriptstyle{CZ}}}}

\def    \Cc    {\operatorname{c^{\circ}_{HZ}}}
\def    \CW    {\operatorname{\bar{c}_{HZ}}}

\def    \chom  {\operatorname{c_{hom}}}
\def    \CHZ   {\operatorname{c_{HZ}}}

\def   \fW    { f_{\scriptscriptstyle{W}} }

\begin{document}


\setlength{\smallskipamount}{6pt}
\setlength{\medskipamount}{10pt}
\setlength{\bigskipamount}{16pt}





\title[Applications of Non-coisotropic Displacement]{Totally non-coisotropic
displacement and its applications to Hamiltonian dynamics}

\author[Ba\c sak G\"urel]{Ba\c sak Z. G\"urel}

\address{Centre de recherches math\'ematiques, 
Universit\'e de Montr\'eal,
P.O. Box 6128, Centre--ville Station,
Montr\'eal, QC, H3C 3J7, CANADA}
\email{gurel@crm.umontreal.ca}

\subjclass[2000]{53D40, 37J45} 

\keywords{Hamiltonian flows, Periodic orbits, Floer homology,
Displacement Principle, Wide Manifold, Action Selector.}

\date{\today}

\bigskip

\begin{abstract}

In this paper we prove the Conley conjecture and the almost existence
theorem in a neighborhood of a closed nowhere coisotropic submanifold
under certain natural assumptions on the ambient symplectic
manifold. Essential to the proofs is a displacement principle for such
submanifolds. Namely, we show that a topologically displaceable
nowhere coisotropic submanifold is also displaceable by a Hamiltonian
diffeomorphism, partially extending the well-known non-Lagrangian
displacement property.

\end{abstract}

\maketitle


\section{Introduction and Main Results}
\labell{main}

In this paper we study Hamiltonian dynamics in a neighborhood of a
nowhere coisotropic submanifold. As a starting point, we establish the
following displacement principle: a closed nowhere coisotropic
submanifold of a symplectic manifold is (infinitesimally) displaceable
provided that there are no topological obstructions to
displaceability.  (In general, a compact subset $M$ of a symplectic
manifold is said to be displaceable if it can be disjoined from itself
by a Hamiltonian diffeomorphism $\varphi_H$, i.e., $\varphi_H(M) \cap
M = \emptyset$. Thus, by definition, a displaceable set is
topologically displaceable.) Then we develop a new version of the
theory of action selectors and use it, combined with the displacement
principle, to prove two results in Hamiltonian dynamics.  Namely, we
prove the Conley conjecture (for non-negative Hamiltonians) and the
almost existence theorem in a neighborhood of a closed nowhere
coisotropic submanifold, under certain natural assumptions on the
ambient manifold. The Conley conjecture, \cite{conley, SZ}, concerns
time-dependent Hamiltonian systems and asserts the existence of
infinitely many periodic points for a Hamiltonian diffeomorphism. The
almost existence theorem due to Hofer and Zehnder and to Struwe,
\cite{hz:ae, ho-ze:capacity, St}, asserts that almost all regular
level sets of a proper autonomous Hamiltonian on $\R^{2n}$ carry
periodic orbits. A similar result has also been proved for $\CP^n$,
symplectic vector bundles, subcritical Stein manifolds, and certain
other symplectic manifolds; see, e.g., \cite{FS, gg3, hv, ke3, lu2,
sc} and also the survey \cite{gi:weinstein} and references therein.
Here, similarly to
\cite{cgk,FS,gi:survey,gg3,gk1,gk2,ke1,ke3,lu,mac1,pol2,felix}, we
focus on these theorems for Hamiltonians supported in a neighborhood
of a closed submanifold.

We also introduce the notion of a \emph{wide} symplectic manifold,
which means that the manifold is open and admits an arbitrarily large
compactly supported Hamiltonian without contractible fast periodic
orbits. Immediate examples of such manifolds include $\R^{2n}$,
cotangent bundles, Stein manifolds and twisted cotangent bundles. The
essence of this property lies in the fact that on a wide manifold the
top degree Floer homology is non-zero for any non-negative compactly
supported Hamiltonian which is not identically zero. This allows us to
construct an action selector for geometrically bounded wide manifolds
and, thus, prove a version of the Conley conjecture.

From now on $W$ will always stand for a wide manifold while $P$ will
denote a general symplectic manifold.

Let us now state the main results of the paper.

\subsection{Displacement Principle}
Let $M$ be a closed connected submanifold of a symplectic manifold
$(P^{2n},\omega)$. We say that $M$ is \emph{nowhere coisotropic} if
$T_xM$ is not a coisotropic subspace of $T_xP$ for any $x \in M$. For
example, a symplectic submanifold is nowhere coisotropic; a
submanifold of middle dimension is nowhere coisotropic if and only if
$\omega|_M$ does not vanish at any point, i.e., $T_xM$ is not a
Lagrangian subspace for any $x \in M $.

Our first result is the following principle which extends or
complements the works of Laudenbach and Sikorav, \cite{LauSi}, and of
Polterovich, \cite{pol1}, and plays a crucial role in the proofs of
the Conley conjecture and the almost existence theorem near nowhere
coisotropic submanifolds.

\begin{Theorem}[Displacement Principle]
\labell{thm:displacement} 
Let $M$ be a closed, connected submanifold of a symplectic manifold
$(P,\omega)$.  Assume that $M$ is nowhere coisotropic and the normal
bundle to $M$ admits a non-vanishing section.  Then $M$ is
infinitesimally displaceable, i.e., there exists a non-vanishing
Hamiltonian vector field which is nowhere tangent to $M$.
\end{Theorem}

When $M$ is of middle dimension, Theorem \ref{thm:displacement} was
proved by Laudenbach and Sikorav, \cite{LauSi}, under a less
restrictive assumption that $M$ is non-Lagrangian and (under the extra
assumption that $TM$ has a Lagrangian complement) by Polterovich,
\cite{pol1}. It has also been known for a long time that $M$ is always
displaceable when $\dim M < n$. (Note that such a submanifold is
automatically nowhere coisotropic.) Thus, Theorem
\ref{thm:displacement} can be thought of as an extension of the
displacement principle to submanifolds of dimension greater than $n$.

In contrast with the middle-dimensional case, the condition that $M$
is nowhere coisotropic cannot be replaced by the requirement that $M$
is (somewhere) non-coisotropic. For a non-coisotropic submanifold can
contain a Lagrangian submanifold and, in this case, $M$ is not
displaceable due to the Lagrangian intersection property. However,
this assumption can possibly be relaxed. We will examine
generalizations of the displacement principle elsewhere.

Theorem \ref{thm:displacement} is proved in Section
\ref{pf:displacement}. Let us now proceed with the applications.

\subsection{The Conley Conjecture}
In its original form, the Conley conjecture asserts that every
Hamiltonian diffeomorphism on $\T ^{2n}$ has infinitely many simple
periodic points, \cite{conley, SZ}. Here ``simple'' means that the
orbits sought are not iterated. A similar conjecture makes sense and
is interesting for other symplectic manifolds too. (Observe that the
example of an irrational rotation on $S^2$ demonstrates that the
conjecture, as stated, fails for manifolds admitting spheres. However,
the statement can be suitably modified to be meaningful and
non-trivial for such manifolds as well; cf. \cite{FrHa}.)

Recently, Ginzburg, \cite{gi:conley}, proved the Conley conjecture for
all closed symplectically aspherical manifolds. Prior to Ginzburg's
work, some particular cases of this conjecture were established. When
the manifold is closed, the conjecture was proved by Salamon and
Zehnder, \cite{SZ}, under the additional assumption that the fixed
points are weakly non-degenerate, and Hingston, \cite{hingston},
established the conjecture for $\T ^{2n}$ without the non-degeneracy
assumption.

Other partial results, not necessarily for closed manifolds, were
obtained under assumptions on the size of the Hamiltonian. Namely, the
conjecture is also known to hold is when the support of the
(time-dependent) Hamiltonian is displaceable. For instance, in $\R
^{2n}$ every compactly supported Hamiltonian has displaceable support,
in which case the conjecture has been proved by Viterbo, \cite{vi2},
and by Hofer and Zehnder, \cite{ho-ze:capacity,hz:book}. Admittedly
this is a very restrictive assumption, especially for closed
manifolds. Yet, this is essentially the only situation in which the
conjecture is known to hold for open manifolds. Under the assumption
that the support is displaceable, the conjecture was proved by
Schwarz, \cite{sc}, for closed symplectically aspherical manifolds and
by Frauenfelder and Schlenk, \cite{FS}, for manifolds that are convex
at infinity. (Recall that a manifold is called convex at infinity if
it is isomorphic at infinity to the symplectization of a compact
contact manifold.) The question is still open for many ``natural''
symplectic manifolds such as cotangent bundles.

We establish this conjecture for Hamiltonians supported in a
neighborhood of a nowhere coisotropic submanifold under certain
assumptions on the ambient manifold. More precisely, we prove

\begin{Theorem}
\labell{thm:cc} 

Assume that $M$ is a closed, nowhere coisotropic submanifold of a
symplectically aspherical manifold $(W,\omega)$. Let $H$ be a non-zero
time-dependent Hamiltonian, supported in a sufficiently small
neighborhood of $M$. 

\begin{itemize}

\item If $W$ is geometrically bounded and wide and $H \geq 0$, then
$H$ has simple (contractible) periodic orbits with positive action and
arbitrarily large period, provided that the time-one map $\varphi_H$
has isolated fixed points with positive action.

\item If $W$ is closed, then $H$ has simple (contractible) periodic
orbits with non-zero action and arbitrarily large period, provided
that the time-one map $\varphi_H$ has isolated fixed points with
non-zero action.

\end{itemize}
\end{Theorem}

We say that $W$ is \emph{wide} if the manifold admits arbitrarily
large, compactly supported, autonomous Hamiltonians $F$ such that all
non-trivial contractible periodic orbits of $F$ have periods greater
than one; see Section \ref{subsec-wide} for a discussion of this
concept. This condition is satisfied for all examples of open
geometrically bounded manifolds known to us.

When $W$ is closed or convex, this theorem can be proved using the
displacement principle, Theorem \ref{thm:displacement}, and the action
selectors introduced in \cite{sc} or \cite{FS} respectively.
Moreover, in these cases it suffices to assume that $W$ is
weakly-exact rather than symplectically aspherical.  However,
the constructions of the action selectors for closed or convex manifolds
do not extend to open manifolds which are merely geometrically
bounded. Hence, the proof of Theorem \ref{thm:cc} for geometrically
bounded manifolds requires developing a new version of the theory of
action selectors; see Section \ref{subsec-actionselector}.

An immediate consequence of Theorem \ref{thm:cc} is the following 
corollary.

\begin{Corollary} 
\labell{cor:cc} 

Let $(W,\omega)$ be geometrically bounded, symplectically aspherical
and wide, and let $M$ be a closed and nowhere coisotropic submanifold
of $W$.  Then, for every non-zero time-dependent Hamiltonian $H \geq
0$ supported in a sufficiently small neighborhood of $M$, the time-one
map $\varphi_H$ has infinitely many simple periodic points
corresponding to contractible periodic orbits of $H$ with positive
action. (A similar statement for closed manifolds has been proved in
\cite{sc}.)

\end{Corollary}

\begin{Remark}
If $W$ is assumed to be convex, the condition that $H \geq 0$ can be
removed both in Theorem \ref{thm:cc} and in Corollary \ref{cor:cc},
provided that the orbits are required only to have non-zero, rather
than positive, action; see \cite{FS}. Moreover, in Theorem
\ref{thm:cc}, the assumption that $W$ is geometrically bounded and
wide can be replaced by the assumption that $W$ admits an exhaustion
$W_1 \subset W_2 \subset \ldots$ by open sets such that each $W_k$ is
symplectomorphic to an open subset of a geometrically bounded and wide
manifold, perhaps depending on $k$.
\end{Remark}

\subsection{The Almost Existence Theorem}
Combining the displacement principle and the results from
\cite{felix}, we prove the following almost existence theorem for
periodic orbits in a neighborhood of a closed nowhere coisotropic
submanifold; see Section \ref{pf:ae}.

\begin{Theorem}
\labell{thm:ae}

Assume that $M$ is a closed, nowhere coisotropic submanifold of a
symplectic manifold $(P,\omega)$ which is geometrically bounded and
strongly semi-positive.  Then the almost existence theorem holds near
$M$: there exists a sufficiently small neighborhood $U$ of $M$ in $P$
such that for any smooth proper Hamiltonian $H \colon U \to \R$, the
level sets $H^{-1} (c)$ carry contractible-in-$P$ periodic orbits of
the Hamiltonian flow of $H$ for almost all $c$ in the range of $H$.

\end{Theorem}

Here $(P^{2n}, \omega)$ is said to be strongly semi-positive if
$c_1(A)\geq0$ for every $A \in \pi_2(P)$ such that $\omega(A)>0$ and
$c_1(A) \geq 2-n$. The condition that $P$ is geometrically bounded
(e.g. convex) is a way to have sufficient control of the geometry of
$P$ at infinity; see Section \ref{prelim} for the definition and
examples.

\begin{Remark}
The displacement results of \cite{felix} rely heavily on \cite{LalMc1,
McDSl}. In Section \ref{proof} we will give a simple proof of this
theorem for symplectically aspherical manifolds $(P,\omega)$ which are
either closed or geometrically bounded and wide.
\end{Remark}

As a particular case, Theorem \ref{thm:ae} implies the almost
existence of periodic orbits in a neighborhood of a closed symplectic
submanifold, provided that $P$ is strongly semi-positive and
geometrically bounded. Note in this connection that almost existence
in a neighborhood of a symplectic submanifold satisfying certain
additional hypotheses was proved by Kerman, \cite{ke3}. On the other
hand, almost existence in a neighborhood of a non-Lagrangian
submanifold of middle dimension was established by Schlenk,
\cite{felix}. Kerman's theorem holds when the ambient manifold $P$ is
symplectically aspherical while Schlenk's requires $P$ to be only
strongly semi-positive. Furthermore, G. Lu, \cite{lu2}, has proved the
almost existence theorem for neighborhoods of symplectic submanifolds
in any symplectic manifold by showing that the contractible
Hofer-Zehnder capacity of such a neighborhood is finite using a deep
and difficult result due to Liu and Tian, \cite{LT2}.

The almost existence theorem is closely related to the existence
problem for periodic orbits of a charged particle in a magnetic field,
also known as the magnetic problem, and to the generalized
Moser-Weinstein theorem; see
\cite{gi:survey,gg3,gk1,gk2,ke1,Mo:orbits}. To be more precise, let
$M$ be a closed Riemannian manifold and let $\eta$ be a closed
two-form (magnetic field) on $M$. Equip $T^*M$ with the twisted
symplectic structure $\omega=\omega_0+\pi^*\eta$, where $\omega_0$ is
the standard symplectic form on $T^*M$ and $\pi\colon T^*M\to M$ is
the natural projection. It is known that $(T^*M,\,\omega)$, a
\emph{twisted cotangent bundle}, is geometrically bounded for any
$\eta$; see \cite{al,cgk,lu}. Finally, let $H$ be the standard kinetic
energy Hamiltonian on $T^*M$. The Hamiltonian flow of $H$ on $W$,
called a \emph{twisted geodesic flow}, is of interest because it
describes, for example, the motion of a charge on $M$ in the magnetic
field $\eta$.

In this setting, as a particular case of Theorem \ref{thm:ae}, we
obtain the existence of contractible twisted geodesics on almost all
low energy levels, provided that the magnetic field is nowhere zero --
a result complementing numerous other theorems on the existence of
twisted geodesics; see, e.g.,
\cite{cgk,gi:survey,gg3,gk1,gk2,ke1,ke3,lu,mac1,pol2,felix}. Note that
the assumption that $\eta$ is nowhere zero ensures that $M$ is nowhere
coisotropic.

\subsection{Organization of the paper} 
In Section \ref{prelim} we set the conventions and notation and recall
relevant results concerning filtered Floer homology and homotopy
maps. The goal of Section \ref{action} is two-fold. We first introduce
and discuss the notion of a wide symplectic manifold. Then we
construct an action selector for wide manifolds which are
geometrically bounded and symplectically aspherical. Here we also
state and prove the properties of this selector. In Section
\ref{proof} we prove the main results of this paper.

\subsection*{Acknowledgments.} The author is deeply grateful to Viktor
Ginzburg for many useful discussions and his numerous valuable remarks
and suggestions. The author also thanks to Yasha Eliashberg, Sam Lisi,
Dusa McDuff, Leonid Polterovich, Tony Rieser, and Aleksey Zinger for
helpful discussions and their suggestions.  Most of this work has been
completed during author's stay at Stony Brook Math Department; she
would like to thank the department for its warm hospitality.


\section{Preliminaries} 
\labell{prelim} 
In this section we set up our conventions and notation and recall the
definition of Floer homology. Here we also define the filtered Floer
homology and examine its dependence on the homotopy of Hamiltonians to
the extent needed in this paper.

We will assume that the manifold $P$ is open, for our proofs will
exclusively focus on the case of open manifolds.

\subsection{Floer homology}
Let $(P,\omega)$ be an open symplectic manifold. In order for the
Floer homology to be defined, we need to impose some additional
conditions on the manifold. To this end, we will always assume that
$P$ is \emph{geometrically bounded}. This assumption gives us
sufficient control of the geometry of $P$ at infinity which is
necessary in the case of open manifolds.  Examples of such manifolds
include symplectic manifolds that are convex at infinity (e.g. compact
symplectic manifolds, $\R^{2n}$, cotangent bundles) as well as twisted
cotangent bundles, which, in general, fail to be convex at
infinity. For the sake of completeness we recall the definition.

\begin{Definition}
\labell{def:gb} 
A symplectic manifold $(P,\omega)$ is said to be \emph{geometrically
bounded} if $P$ admits an almost complex structure $J$ and a complete
Riemannian metric $g$ such that
\begin{itemize}
\item $J$ is uniformly $\omega$-tame, i.e.,
for some positive constants $c_1$ and $c_2$ we have
$$
\omega(X, JX)\geq c_1\| X\|^2 \quad\text{and}\quad
|\omega(X, Y)|\leq c_2\| X\|\,\|Y\|
$$
for all tangent vectors $X$ and $Y$ to $P$;
\item the sectional curvature of $(P,g)$ is bounded from above and
the injectivity radius of $(P,g)$ is bounded away from zero.
\end{itemize}
\end{Definition}

We refer the reader to \cite{al, cgk, lu} for a discussion of
geometrically bounded manifolds.  In particular, though we have not
yet recalled the definition of Floer homology, let us note that the
compactness theorem for the moduli spaces of Floer's trajectories for
open geometrically bounded manifolds holds; this is a consequence of
Sikorav's version of the Gromov compactness theorem; see \cite{al}.

Furthermore, assume that $(P,\omega)$ is \emph{symplectically
aspherical}, i.e.,
$$
\omega|_{\pi_2(P)}=0 \quad\text{and}\quad c_1(TP)|_{\pi_2(P)}=0.
$$ 
We will indicate when $P$ need not be symplectically aspherical, as is
the case in Theorem \ref{thm:ae}.

Among manifolds which are symplectically aspherical and geometrically
bounded are $\R^{2n}$, symplectic tori, cotangent bundles and twisted
cotangent bundles when the form on the base is weakly exact. Under
these hypotheses, the filtered $\Z$-graded Floer homology of a
compactly supported Hamiltonian on $P$ is defined as follows.

Recall that for a time-dependent Hamiltonian $H\colon S^1 \times P\to
\R$, the action functional on the space of smooth contractible loops
$\Lambda P$ is defined as
\begin{equation} 
\labell{eq:action}
A_H(x) = - \int_{D^2} \bar{x}^* \omega +
\int_{S^1}H(t,x)\,dt,
\end{equation}
where $x\colon S^1 \to P$ is a contractible loop and $\bar{x}\colon
D^2\to P$ is a map of a disk, bounded by $x$, and $S^1 = \R /
\Z$. Since $P$ is symplectically aspherical $A_H(x)$ is
well-defined. The Hamiltonian vector field $X_H$ is defined by the
equation $ i_{X_H} \omega = -dH $. Let $\varphi^t_H$ denote the
time-dependent flow of $X_H$ and, in particular,
$\varphi_H=\varphi^1_H$ denote the time-one flow.

By the least action principle, the critical points of $A_H$ are
exactly contractible one-periodic orbits of the Hamiltonian flow of
$H$. We denote by $\Pp_H$ the collection of such orbits and let
$\Pp_H^{(a,\,b)} \subset \Pp_H$ stand for the collection of orbits
with action in the interval $(a,\,b)$. The action spectrum $\Ss(H)$ of
$H$ is the set of critical values of $A_H$. In other words, $\Ss(H)=\{
A_H(x)\mid x \in \Pp_H\}$. This is a zero measure set; see, e.g.,
\cite{hz:book,sc}.

Throughout this paper we will assume that $H$ is compactly supported
and set $\supp H=\bigcup_{t\in S^1}\supp H_t$. In this case, $\Ss(H)$
is closed and hence nowhere dense.

Let $J=J_t$ be a time-dependent almost complex structure on $P$. A Floer
anti-gradient trajectory $u$ is a map $u\colon \R\times S^1\to P$
satisfying the equation
\begin{equation}
\label{eq:floer}
\frac{\p u}{\p s}+ J_t(u) \frac{\p u}{\p t}=-\nabla H_t(u).
\end{equation}
Here the gradient is taken with respect to the time-dependent
Riemannian metric $\omega(\cdot,J_t\cdot)$.  Denote by $u(s)$ the
curve $u(s,\cdot)\in \Lambda P$.

The energy of $u$ is defined as
\begin{equation}
\label{eq:energy}
E(u)=\int_{-\infty}^\infty \left\|\frac{\p u}{\p s}\right\|_{L^2(S^1)}^2\,ds
=\int_{-\infty}^\infty \int_{S^1}\left\|\frac{\p u}{\p t}-J\nabla H (u)
\right\|^2 \,dt\,ds.
\end{equation}
We say that $u$ is asymptotic to $x^\pm\in \Pp_H$ as $s\to\pm \infty$,
or connecting $x^-$ and $x^+$, if $\lim_{s\to\pm\infty} u(s)=x^\pm$ in
$\Lambda P$. In this case
$$
A_H(x^-)-A_H(x^+)=E(u).
$$
We denote the space of Floer trajectories connecting $x^-$ and $x^+$,
with the topology of uniform $C^\infty$-convergence on compact sets,
by $\Mm_H(x^-,x^+,J)$. This space carries a natural $\R$-action
$(\tau\cdot u)(t,s)=u(t,s+\tau)$ and we denote by
$\hat{\Mm}_H(x^-,x^+,J)$ the quotient $\Mm_H(x^-,x^+,J)/\R$.

Recall that $x \in \Pp_H$ is said to be non-degenerate if
$d\varphi_H\colon T_{x(0)}P\to T_{x(0)}P$ does not have one as an
eigenvalue. In this case, the so-called Conley--Zehnder index
$\MUCZ(x)\in\Z$ is defined; see, e.g., \cite{sa,SZ}.  Here we
normalize $\MUCZ$ so that $\MUCZ(x)=n$ when $x$ is a non-degenerate
maximum of an autonomous Hamiltonian with a small Hessian. Assume that
all periodic orbits with actions in the interval
$[A_H(x^+),\,A_H(x^-)]$, including $x^\pm$, are non-degenerate.  Then,
for a generic $J$, suitable transversality conditions are satisfied
and $\Mm_H(x^-,x^+,J)$ is a smooth manifold of dimension
$\MUCZ(x^-)-\MUCZ(x^+)$; see, e.g., \cite{fh,SZ} and references
therein.

\subsection{Filtered Floer homology} 
In this section we briefly outline the construction of filtered Floer
homology following closely \cite{gi:coiso}; see also
\cite{bps,cgk,fh,gg3,sc}.

Throughout the discussion of the filtered Floer homology
$\HF^{(a,\,b)}(P)$, we assume that all intervals are in the positive
range of actions, i.e., $a>0$ for any interval $(a,\,b)$. This
condition is clearly necessary, for $H$ is assumed to be compactly
supported and thus it always has trivial degenerate periodic orbits if
$P$ is open.

\subsubsection{Filtered Floer homology: definitions}
Let $H$ be a compactly supported Hamiltonian on an open symplectic
manifold $P$ which is symplectically aspherical and geometrically
bounded. Assume that all contractible one-periodic orbits of $H$ with
positive action are non-degenerate. This is a generic condition.
Consider an interval $(a,\,b)$, with $a>0$, such that $a$ and $b$ are
outside $\Ss(H)$. Then the collection $\Pp^{(a,\,b)}_H$ is finite.
Assume furthermore that $J$ is regular, i.e., the necessary
transversality conditions are satisfied for moduli spaces of Floer
trajectories connecting orbits from $\Pp^{(a,\,b)}_H$. This is again a
generic property as can be readily seen by applying the argument from
\cite{fh,FHS,SZ}.

Let $\CF_k^{(a,\,b)}(H)$ be the vector space over $\Z_2$ generated by
$x\in \Pp^{(a,\,b)}_H$ with $\MUCZ(x)=k$. Define
$$
\p \colon \CF_k^{(a,\,b)}(H)\to \CF_{k-1}^{(a,\,b)}(H)
$$
by
$$
\p x=\sum_y \#\big(\hat{\Mm}_H(x,y,J)\big)\cdot y,
$$
where the summation extends over all $y\in \Pp^{(a,\,b)}_H$ with
$\MUCZ(y)=\MUCZ(x)-1$ and $\#\big(\hat{\Mm}_H(x,y,J)\big)$ is the
number of points, modulo 2, in $\hat{\Mm}_H(x,y,J)$. (Recall that in
this case $\hat{\Mm}_H(x,y,J)$ is a finite set by the compactness
theorem.) Then, as is well known, $\p^2=0$.  The resulting complex
$\CF^{(a,\,b)}(H)$ is the filtered Floer complex for $(a,\,b)$. Its
homology $\HF^{(a,\,b)}(H)$ is called the filtered Floer homology.
This is essentially the standard definition of Floer homology with
critical points having action outside $(a,\,b)$ being ignored.  In
general, $\HF^{(a,\,b)}(H)$ depends on the Hamiltonian $H$, but not on
$J$; see, e.g., \cite{gi:coiso}.

\begin{Remark}
It is clear that the same construction, with suitable modifications,
works for closed manifolds. In this case,
$\HF(H)=\HF^{(-\infty,\infty)}(H)$ is the ordinary Floer homology.
Moreover, $\HF_*(H)=H_{*+n}(P;\,\Z_2)$.
\end{Remark}

Let $a<b<c$. Assume that all of the above assumptions are satisfied for
all three intervals $(a,\,c)$ and $(a,\,b)$ and $(b,\,c)$. Then clearly
$\CF^{(a,\,b)}(H)$ is a subcomplex of $\CF^{(a,\,c)}(H)$, and
$\CF^{(b,\,c)}(H)$ is naturally isomorphic to the quotient complex
$\CF^{(a,\,c)}(H)/\CF^{(a,\,b)}(H)$. As a result, we have the long exact
sequence
\begin{equation}
\label{eq:seq}
\ldots\to \HF^{(a,\,b)}(H)\to\HF^{(a,\,c)}(H)\to \HF^{(b,\,c)}(H)\to\ldots.
\end{equation}

In the construction of the action selector for open manifolds given in
Section~\ref{sel-defn}, we will work with filtered Floer homology for
the interval $(0,\,b)$ even though $0$ is necessarily a critical value
of the action functional. This homology is defined as
\begin{equation}
\label{eq:zero-infty}
\HF^{(0,\,b)}(H)=\varprojlim_{\eps\to 0+}\HF^{(\eps,\,b)}(H),
\end{equation}
where the inverse limit is taken with respect to the quotient maps and
$\eps\to 0+$ in the complement of $\Ss(H)$.  (It is clear that this
definition is equivalent to the original one when $P$ is closed and
$0$ is not in $\Ss(H)$.)

\subsubsection{Homotopy}
\label{sec:homotopy}
Let us now examine the dependence of $\HF^{(a,\,b)}(H)$ on $H$.
Consider a homotopy $H^s$ of Hamiltonians from $H^0$ to $H^1$. By
definition, this is a family of Hamiltonians parametrized by $s\in
\R$, and such that $H^s\equiv H^0$ when $s$ is large negative and
$H^s\equiv H^1$ when $s$ is large positive. Furthermore, let $J^s$ be
a family of $t$-dependent almost complex structures such that again
$J^s\equiv J^0$ when $s\ll 0$ and $J^s\equiv J^1$ when $s\gg 0$. (We
will most of the time suppress $J^s$ in the notation and refer to
$H^s$ as the homotopy.)

For $x\in \Pp^{(a_0,\,b_0)}_{H^0}$ and $y\in
\Pp^{(a_1,\,b_1)}_{H^1}$ denote by $\Mm_{H^s}(x,y,J^s)$ the space of
solutions of \eqref{eq:floer} with $H=H^s$ and $J=J^s$.

The regularity property takes the following form for open manifolds:
$(H^s,J^s)$ is said to be regular if the transversality requirements
are met along all homotopy trajectories connecting periodic orbits
with positive action. This is a generic property as can be seen by
arguing as in \cite{fh,FHS,SZ}. (When $P$ is closed, regularity of a
homotopy $(H^s,J^s)$ is understood in the standard sense, i.e., the
standard transversality requirements are met by the homotopy
$(H^s,J^s)$; see \cite{fh,FHS,SZ}.)

When the transversality conditions are satisfied, $\Mm_{H^s}(x,y,J^s)$
is a smooth manifold of dimension $\MUCZ(x)- \MUCZ(y)$. In particular,
$\Mm_{H^s}(x,y,J^s)$ is a finite set when $\MUCZ(x)=\MUCZ(y)$. Define
the homotopy map
$$
\Psi_{H^0H^1}\colon \CF^{(a_0,\,b_0)}(H^0)\to \CF^{(a_1,\,b_1)}(H^1)
$$
by
$$
\Psi_{H^0H^1}( x)=\sum_y \#\big(\Mm_{H^s}(x,y,J^s)\big)\cdot y.
$$
Here the summation is over all orbits $y\in \Pp^{(a_1,\,b_1)}_{H^1}$ with
$\MUCZ(y)=\MUCZ(x)$ and $\#\big(\Mm_{H^s}(x,y,J^s)\big)$
is the number of points, modulo 2, in this moduli space.

The map $\Psi_{H^0H^1}$ depends on the entire homotopy $(H^s,J^s)$ and
in general is \emph{not} a map of complexes. However, $\Psi_{H^0H^1}$
becomes a homomorphism of complexes when $(a_0,\,b_0) = (a_1,\,b_1)$
and the homotopy is monotone decreasing, i.e., $\p_s H^s\leq 0$
point-wise. Moreover, the induced map in homology is then independent
of the homotopy, within the class of decreasing homotopies, and
commutes with the maps from the exact sequence \eqref{eq:seq}.  (The
reader is referred to, e.g., \cite{bps,cgk,fh,sa,SZ,sc,vi:functors},
for the proofs of these facts for both open and closed manifolds.)
There are other instances when the same is true. In particular, this
is the case when the location of the intervals $(a_0,\,b_0)$ and
$(a_1,\,b_1)$ is compatible with the growth of the Hamiltonians in the
homotopy, as streamlined by the following theorem from
\cite{gi:coiso}; see also \cite{sc}. (This theorem holds for both open
and closed manifolds.)

\begin{Theorem}[\cite{gi:coiso}]
\labell{thm:homotopy}
Let $H^s$ be a homotopy such that 
$$
\int_{-\infty}^{\infty}\int_{S^1}\max_{P} \p_s H^s_t \,dt\,ds \leq C,
$$
where $C \in \R$.  Then
$$
\Psi_{H^0H^1}\colon \CF^{(a,\,b)}(H^0)\to \CF^{(a+C,\,b+C)}(H^1)
$$
is a homomorphism of complexes for any interval $(a,\,b)$.  Hence,
$\Psi_{H^0H^1}$ induces a map in Floer homology, also denoted by
$\Psi_{H^0H^1}$. This map sends the exact sequence \eqref{eq:seq} for
$H^0$ and $(a,\,b,\,c)$ to the exact sequence \eqref{eq:seq} for $H^1$
and $(a+C,\,b+C,\,c+C)$, i.e., on the level of homology
$\Psi_{H^0H^1}$ commutes with all maps in the long exact sequence
\eqref{eq:seq}.

\end{Theorem}


\section{Action selectors in wide manifolds}
\labell{action}

The proof of Theorem \ref{thm:cc} relies on the theory of action
selectors, which is one of the standard approaches to the problem,
\cite{FS, hz:book, sc, vi2}. This theory is well developed for
weakly-exact closed or convex manifolds; see, e.g., \cite{FGS,FS,sc} and
also \cite{oh3} for the theory in a more general setting. The main
ingredient in these constructions of an action selector is an
identification between Floer homology spaces for different
Hamiltonians.  For example, for symplectically aspherical closed or
convex manifolds, the Floer homology of any Hamiltonian is isomorphic
to the homology of the manifold. Then an action selector can be
associated to any homology class.  However, these constructions do not
generalize to the case of open manifolds which are merely
geometrically bounded. The main obstacle is that the Floer homology
for Hamiltonians on such manifolds is no longer an invariant, i.e., it
depends on the Hamiltonian.  For the homology can be defined only for
action intervals that do not contain zero and hence, in contrast with
the case of closed or convex manifolds, there is, in general, no
relation between the Floer homology for different
Hamiltonians. Nevertheless, to construct an action selector it
suffices to have a class in the Floer homology of a Hamiltonian, which
is in some sense canonical, even though the homology group it belongs
to depends on the Hamiltonian. In this section we construct an action
selector for geometrically bounded wide manifolds. The homology class
for which the selector is defined is essentially the fundamental class
of the manifold modulo infinity. Then this selector, for non-negative
Hamiltonians, has properties similar to those of action selectors
constructed in \cite{FS,sc}. 

\subsection{Wide manifolds}
\labell{subsec-wide} 
Let us now introduce and discuss the notion of a \emph{wide}
symplectic manifold.

\begin{Definition}
\labell{def:wide1} 
A symplectic manifold $(W,\omega)$ is said to be wide if it is open
and if for every constant $C \geq 0$ and for every compact subset $X
\subset W$, there exists a function $K \colon W \to [0,\infty)$ such
that

\begin{itemize}
\item[(W1)] $K$ is compactly supported; 
\item[(W2)] $K|_X \geq C$; 
\item[(W3)] the Hamiltonian flow of $K$ has no non-trivial contractible
fast periodic orbits.  
\end{itemize}
\end{Definition}

We call a non-trivial orbit fast if its period is less than or equal
to one. Otherwise an orbit will be called slow. 

\begin{Remark}
The condition (W2) can be replaced by the condition $\mathrm {(W2')}
\colon \max \,K \equiv K|_X \equiv C$.  For, as explained in
\cite{gg3, hz:book}, one can cut off a function meeting the
requirement (W2) without creating new fast periodic orbits and produce
another function satisfying the condition $\mathrm{(W2')}$.
\end{Remark}

Examples of wide manifolds include cotangent bundles and Stein
manifolds, or more generally symplectic manifolds that are convex at
infinity. More importantly, twisted cotangent bundles, which are
geometrically bounded but, in general, fail to be convex at infinity,
are wide. Non-compact covering spaces of compact manifolds are also
examples of wide manifolds, \cite{gi:private}.  Furthermore, the
product of two wide manifolds is wide and so is the product of a
compact and a wide manifold. On the other hand, a manifold with finite
contractible Hofer-Zehnder capacity cannot be wide, as easily follows
from Definition \ref{def:wide1}; see, e.g., Section \ref{pf:ae} for
the definition of this capacity.

\begin{Remark}
Wideness, although indispensable for our construction, appears to be a
rather mild assumption: the author is not aware of any example of an
open geometrically bounded manifold which is not wide. It is not yet
clear, however, whether every open geometrically bounded manifold is
wide. The relation between geometrical boundedness and the concept of
wideness is interesting but quite complicated, and we will address
this question elsewhere.
\end{Remark} 

The property of wideness can also be viewed in terms of the restricted
relative Hofer-Zehnder capacity introduced in \cite{gg3} or as the
property of a manifold to admit a slow proper function. Namely, we
have

\begin{Proposition}
\labell{prop:equivalent}
Let $(W,\,\omega)$ be a symplectic manifold. Then the following
statements are equivalent.
\begin{itemize}
\item[(i)] $(W,\,\omega)$ is wide.
\item[(ii)] $\CW(W,X)= \infty$ for every compact subset $X \subset W$,
where $\CW(W,X)$ is the restricted (or contractible) relative
Hofer-Zehnder capacity introduced in \cite{gg3}. 
\item[(iii)] $(W,\,\omega)$ admits a non-negative proper function
without non-trivial contractible fast periodic orbits.
\end{itemize}
\end{Proposition}

We omit the proof of Proposition \ref{prop:equivalent}, for it is
essentially straightforward and we will mainly be using Definition
\ref{def:wide1}.

\begin{Remark}
One could also replace the condition $\mathrm{(W3)}$ in Definition
\ref{def:wide1} by the one that the Hamiltonian flow of $K$ has no
non-trivial fast periodic orbits, hence dropping the requirement that
the orbits be contractible. However, this would be a more restrictive
requirement on $W$; for instance, a cylinder of finite area, which is
wide according to Definition \ref{def:wide1}, would then fail to be
wide.
\end{Remark}

Let us now turn to constructing an action selector for geometrically
bounded wide manifolds.

\subsection{An action selector for wide manifolds}
\labell{subsec-actionselector}

We assume that $(W^{2n},\omega)$ is a symplectically aspherical,
geometrically bounded and wide manifold. We will construct an action
selector for non-negative compactly supported Hamiltonians on~$W$.

\subsubsection{The definition}
\labell{sel-defn}

Let $H \colon S^1 \times W \to \R$ be a compactly supported
Hamiltonian such that $H \geq 0$.  It is easy to see that, since $W$
is wide, there exists a smooth compactly supported function $F \colon
W \to [0,\,\infty)$ without non-trivial contractible fast periodic
orbits and such that $F \ge H$ point-wise.  (This is essentially the
definition of a wide manifold.)  Without loss of generality we may
assume that $\supp (F)$ is a smooth connected manifold with boundary
and that $F$ is a Morse function with finitely many critical points
when restricted to the interior of the support. (These requirements
are generic.) From now on we will call such functions \emph{wide} and
we will reserve the notation $F$ for them.

Under these assumptions, 
\begin{equation}
\label{eq:floer=morse}
\HF_*^{(0,\,\infty)}(F) \cong \HM_{*+n}^{(0,\,\infty)}(F),
\end{equation}
where $\HM_*^{(0,\,\infty)}(F)$ and $\HF_*^{(0,\,\infty)}(F)$ denote,
respectively, the (filtered) Morse and Floer homology of $F$ for the
interval $(0, \, \infty)$.

For the sake of completeness, let us explain first why this
isomorphism holds. Intuitively, this is obvious. For $F$ has no
non-trivial contractible one-periodic orbits by definition and, hence,
Floer and Morse complexes are the same as vector spaces. Note,
however, that the two differentials may be different. We,
nevertheless, claim that the resulting homology groups are isomorphic.

To this end, let $\epsilon>0$ be small enough so that $\epsilon F$ is
$C^2$-small. Recall that for $C^2$-small functions, when $W$ is
symplectically aspherical, Morse and Floer homology groups are isomorphic. Then
$$
\HF_*^{(0,\,\infty)}(\epsilon F) \cong \HM_{*+n}^{(0,\,\infty)}(\epsilon F).
$$
Consider a monotone decreasing homotopy $F_s$, for $s \in [0,\,1]$,
from $F$ to $\epsilon F$. Choose $\delta >0$ such that it is below any
critical values of $\epsilon F$. As a result, action values for $F_s$
for all $s \in [0,\,1]$ are greater than $\delta > 0$. By the homotopy
invariance of Floer homology, we then have the isomorphism
$$
\HF_*^{(\delta,\,\infty)}(F) \cong \HF_*^{(\delta,\,\infty)}(F_s)
\text{ for all }s\in [0,\,1].
$$
Also note that $\HF_*^{(0,\,\infty)}(F_s) =
\HF_*^{(\delta,\,\infty)}(F_s)$ for all $s \in [0,\,1]$ since zero is an
isolated critical value and $\delta$ is smaller than any critical
value of $F_s$.

Finally, we have
$$
  \HF_*^{(0,\,\infty)}(F)
\cong \HF_*^{(0,\,\infty)}(\epsilon F)
\cong \HM_{*+n}^{(0,\,\infty)}(\epsilon F)
\cong \HM_{*+n}^{(0,\,\infty)}(F).
$$

\begin{Remark}
Observe that $\HM_{*+n}^{(0,\,\infty)}(F)$ is just the ordinary
homology of $\supp (F)$ modulo its boundary, i.e.,
$\HM_{*+n}^{(0,\,\infty)}(F)=H_{*+n}(\supp F, \partial(\supp F);
\Z_2)$.
\end{Remark}

Let us now turn to the definition of the action selector. Since $\supp
(F)$ is connected, equation \eqref{eq:floer=morse}, in particular,
implies that
$$
\HF_n^{(0,\,\infty)}(F) \cong \HM_{2n}^{(0,\,\infty)}(F) \cong {\Z}_2.
$$
Denote by $[\max _F]$ the generator of $\HF_n^{(0,\,\infty)}(F) \cong
\Z_2$, which can be thought of as the fundamental class. Consider the
image of this class under the monotonicity map
$$
\Psi_{\scriptscriptstyle{FH}} \colon 
\HF_n^{(0,\,\infty)}(F) \to \HF_n^{(0,\,\infty)}(H),
$$
induced by a monotone decreasing homotopy from $F$ to $H$. (This map
is independent of the choice of the homotopy, as discussed in Section
\ref{prelim}.) Using the homotopy invariance of the Floer homology it
is not hard to show that $\Psi_{\scriptscriptstyle{FH}}[\max F]$ is
independent of the choice of $F$; denote this class by $[\max]$.

\begin{Definition} 
\labell{defn:sel}
The action selector is defined as
$$\sigma(H)=\mathrm{inf} \, \{a > 0 \,\,|\,\,j_a^H [\max] =0 \},$$
where
$$
j_a^H \colon \HF_n^{(0,\,\infty)}(H) \to \HF_n^{(a,\,\infty)}(H),
$$
is induced by the quotient map.
\end{Definition}

\begin{Remark}
Note that this definition makes sense for arbitrary Hamiltonians,
i.e., we need not assume that $H$ is non-negative. However, it is not
clear whether $\sigma$ would then have the properties (AS0)-(AS6)
listed below, which are crucial for applications. Also, when $W$ is
convex, this selector is equal to the one constructed in \cite{FS}.
\end{Remark}

\subsubsection{Properties of the action selector}
\labell{sel-prop}

Let $\Ss(H)$ denote the action spectrum of $H$ and let $\Ham_c^+(W)$
denote the cone in the group of compactly supported Hamiltonian
diffeomorphisms of $W$ generated by non-negative Hamiltonians. The
action selector $\sigma \colon \Ham_c^+(W) \to [0,\,\infty)$
constructed above has the following properties:

\begin{itemize}

\item[(AS0)] $\sigma(H)$ is a spectral value: $\sigma(H) \in \Ss(H)$;

\item[(AS1)] $\sigma(H)$ is monotone:
if $0 \le H \leq K$ then $0 \le \sigma(H) \leq \sigma(K)$;

\item[(AS2)] $\sigma(H)$ is non-degenerate:
$0< \sigma(H) < \infty$ if $H \not\equiv 0$;

\item[(AS3)] $\sigma(H)$ is continuous in $H$ with respect to the Hofer
norm, and, in particular, it is $C^0$-continuous;

\item[(AS4)] $\sigma(H) \leq \int_0^1 \max {H_t} \,dt = \| H \|$,
where $||\, \cdot \,||$ denotes the Hofer norm;

\item[(AS5)] if $H$ is autonomous and has no non-trivial contractible
fast periodic orbits, then $\sigma(H) = \max H $;

\item[(AS6)] if $H$ and $K$ generate the same time-one flow
and ${\varphi}_H^t$ and ${\varphi}_K^t$ are homotopic (with fixed end
points) in $\Ham_c ^+(W)$, then $\sigma(H)=\sigma(K)$;

\item[(AS7)] $\sigma(H) \leq e_V$ for any $H$ supported in $V\subset
W$, where $e_V$ denotes the displacement energy of $V$. Thus,
$\sigma(H)$ is \emph{a priori} bounded from above by $e_V<\infty$ if
$V$ is displaceable.

\end{itemize}

\subsubsection{Proofs of the properties of the selector}
\labell{proof-sel-prop} 
We will prove these properties in varying degree of detail; some
proofs are very similar to those for the selectors constructed in
\cite{FS,sc}, whereas some proofs require modifications. We will
mainly focus on the new parts and refer to the literature for the
standard ones.

\emph{(AS0)} First recall that $\Ss(H)$ is compact and nowhere dense.

Assume the contrary: $\sigma(H) \not\in \Ss(H)$. Then, since $\Ss(H)$
is compact, for a small enough number $\delta >0$ we have
$[\sigma(H)-\delta, \, \sigma(H)+\delta] \cap \Ss(H) =
\emptyset$. Hence, by the definition of the selector, there exists a
number $a$ with $ \sigma(H) < a \le \sigma(H)+\delta $ such that
${j_a^H}[\max ] =0$.

Let $c$ be such that $\sigma(H)-\delta \le c < \sigma(H) $. Then we
have the isomorphism
$$\HF_n^{(a,\,\infty)}(H) \cong \HF_n^{(c,\,\infty)}(H),$$ since there
are no critical values of $A_H$ in $[c,\,a]$, and the diagram
\begin{displaymath}
\xymatrix
{
& 
  \HF_n^{(0,\,\infty)}(H) \ar[dl]_{0=j_a^H}  \ar[dr]^{j_c^H} &  \\
  \HF_n^{(a,\,\infty)}(H) \ar[rr]^{\cong} &
& \HF_n^{(c,\,\infty)}(H)
}
\end{displaymath}
commutes. Thus ${j_c^H}[\max ] =0$. This contradicts the definition of
$ \sigma(H)$.  We conclude that $\sigma(H) \in \Ss(H)$.

\emph{(AS1)} Let $K \ge H \ge 0$ and let $F$ be a wide function such
that $F \ge K \ge H$. Monotonicity is a consequence of the
commutativity of the following diagram:

\begin{displaymath}
\xymatrix
{
  \HF_n^{(0,\,\infty)}(F) \ar[r]^{\Psi_{FH}} \ar[dr]_{\Psi_{FK}} 
& \ar@{}[dr] \HF_n^{(0,\,\infty)}(H) \ar[r]^{j_a^H}  
& \HF_n^{(a,\,\infty)}(H) \\
& \HF_n^{(0,\,\infty)}(K) \ar[r]_{j_a^K}  
& \HF_n^{(a,\,\infty)}(K) \ar[u]_{\Psi_{KH}}
}
\end{displaymath}

\emph{(AS2)} Finiteness of the selector follows from the compactness
of $\Ss(H)$, for $\HF_n^{(a,\,\infty)}(H) = 0$ for any $a >
\sup \Ss(H)$.

Proving the non-degeneracy, i.e., $\sigma(H)>0$ for any $H \not\equiv
0$, requires more work. First observe that we can find a $C^2$-small
\emph{space-time} bump function $f\not\equiv 0$ such that $0 \le f \le
H$ for all $t\in S^1$. This is simply because $H \geq 0$ and $H
\not\equiv 0$ for some $t$. Hence, by the monotonicity of the selector,
it suffices to show that $\sigma(f)> 0$.

More precisely, let $f(t,x)= f_{\scriptscriptstyle{S^1}}(t) \cdot
\fW(x) \colon S^1 \times W \to [0,\,\infty)$ be a space-time bump
function satisfying $0 \le f \le H$ for all $t\in S^1$. Here
$\fW(x)\colon W \to [0,\,\infty)$ and $f_{\scriptscriptstyle{S^1}}(t)
\colon S^1 \to [0,\,\infty)$ are both bump functions in the usual
sense, and $\fW$ is autonomous and $C^2$-small. The time-one flow of
$f$ differs from the time-one flow of $\fW$ only by a positive factor
equal to the integral of $f_{\scriptscriptstyle{S^1}}$ over the
circle.  Let us assume, for the sake of simplicity, that $ \int_{S^1}
f_{\scriptscriptstyle{S^1}}(t) \, dt =1$.  This can be achieved by
choosing $\fW$ sufficiently small so that $f$ still fits underneath
$H$. Then the action spectra of $f$ and $\fW$ are the same.

We claim that $\HF_*^{(a,\,b)}(f) \cong \HF_*^{(a,\,b)}(\fW)$ for all
positive intervals of action $(a,\,b) \subseteq (0,\,\infty)$ such
that $a$ and $b$ are not in $\Ss(f)=\Ss(\fW)$. To see this, let $K_s =
s \fW + (1-s) f$ be the linear homotopy from $f$ to $\fW$ for $s \in
[0,\,1]$. It is easy to see that all Hamiltonians in this homotopy
have the same time-one flow as that of $\fW $ and, hence, the only
critical points of $A_{K_s}$ are constant one-periodic orbits, i.e.,
the critical points of $\fW$. Thus the action spectrum $\Ss(K_s)$ for
any $s$ consists of two action values: zero and $\max \fW$. This
implies that for $a,\,b \not\in \Ss(K_s)$ no periodic orbit with
action outside the range $(a,\,b)$ will enter or exit this interval
during the course of the homotopy. In this case the Floer homology
groups are isomorphic for all $K_s$; see, e.g., \cite{bps, gi:coiso,
vi:functors}. In particular, $\HF_*^{(a,\,b)}(f) \cong
\HF_*^{(a,\,b)}(\fW)$. (Note that this isomorphism is not induced by
the homotopy map since $K_s$ is not a monotone homotopy.)

We next show that $\sigma(\fW)=\max \fW >0$. To this end, let $F$ be a
wide function satisfying $F \ge \fW$ and recall that
$\HF_*^{(0,\,\infty)} (F) \cong \HM_{*+n}^{(0,\,\infty)}
(F)$. Moreover, since $\fW$ is a $C^2$-small bump function, Morse and
Floer homology groups are isomorphic: $\HF_*^{(0,\,\infty)}(\fW) \cong
\HM_{*+n}^{(0,\,\infty)}(\fW)$. Then the diagram
\begin{displaymath}
\xymatrix 
{
{\Z}_2 \cong \HF_n^{(0,\,\infty)}(F) 
\ar[r]^{\cong}
\ar[d]_{\Psi_{F \fW}}  &
\HM_{2n}^{(0,\,\infty)}(F) \cong {\Z}_2
\ar[d]^{\Psi_{F \fW}} \\
\HF_n^{(0,\,\infty)}(\fW) 
\ar[r]^{\cong} 
\ar[d]_{j_a^{\fW}} &
\HM_{2n}^{(0,\,\infty)}(\fW)
\ar[d]^{j_a^{\fW}} \\
\HF_n^{(a,\,\infty)}(\fW) 
\ar[r]^{\cong} &
\HM_{2n}^{(a,\,\infty)}(\fW)
}
\end{displaymath}
is commutative. This can easily be seen by factoring the horizontal
isomorphisms through isomorphisms induced by monotone homotopies and
using the fact that all such diagrams commute when the functions
involved are $C^2$-small.

Let us focus on the right-hand side of the diagram. By Morse theory,
the map
\begin{displaymath}
\xymatrix
{
{\Z}_2 \cong \HM_{2n}^{(0,\,\infty)}(F)
\ar[r]^{\Psi_{F\fW}} 
&
\HM_{2n}^{(0,\,\infty)}(\fW) \cong {\Z}_2
}
\end{displaymath}
is non-zero, and sends $[\max _F]$ to $[\max _{\fW}]$. Moreover,
$j_a^{\fW}[\max _{\fW}]=[\max \fW]$ for any positive $a < \max \fW$
and $j_a^{\fW}[\max _{\fW}]=0$ for any $a \ge \max \fW$. Commutativity
of the diagram then implies that $\sigma(\fW)=\max \fW >0$.

As the last step observe that $\sigma(f)=\sigma(\fW)=\max \fW > 0$,
which finishes the proof of non-degeneracy.  To see this, note that
the diagram
\begin{displaymath}
\xymatrix
{
  \HF_n^{(0,\,\infty)}(F) \ar[r]^{\Psi_{Ff}} 
\ar[dr]_{\Psi_{F\fW}} 
& \ar@{}[dr]
  \HF_n^{(0,\,\infty)}(f) \ar[r]^{j_a^f}  \ar@{=}[d]
& \HF_n^{(a,\,\infty)}(f) \ar@{=}[d] \\
& \HF_n^{(0,\,\infty)}(\fW) 
\ar[r]_{j_a^{\fW}}  
& \HF_n^{(a,\,\infty)}(\fW) 
}
\end{displaymath}
is also commutative, where $F$ is a wide function satisfying $F\ge f$
and $F \ge \fW$.

\begin{Remark}
An important consequence of non-degeneracy of $\sigma$ is that
$\HF_n^{(0,\,\infty)}(H) \neq 0$ for any non-negative Hamiltonian $H
\not\equiv 0$. To see this, note that the \emph{non-zero} map
$\HF_n^{(0,\,\infty)}(F) \to \HF_n^{(0,\,\infty)}(\fW)$ can be
factored through $\HF_n^{(0,\,\infty)}(H)$, where $F$ and $\fW$ are as
in the proof above. This fact is used in \cite{gi:coiso}.
\end{Remark}

\emph{(AS3)} Recall that this property asserts the continuity of the
selector with respect to the Hofer norm. Namely, we claim that
$$
|\sigma(H)-\sigma(K)| \le \| H-K \| \text{ for non-negative } H \text{
 and } K.
$$

Note first that it suffices to prove continuity for non-degenerate
Hamiltonians.  Since we are assuming that $W$ is open, this means that
all one-periodic orbits with positive action are non-degenerate.

Keeping the notation from Section~\ref{sec:homotopy}, let $K=H^0$ and
$H=H^1$, and consider a linear homotopy $H^s$ from $H^0$ to $H^1$,
i.e., $H^s=(1-\phi(s)) H^0 +\phi(s) H^1$, where $\phi \colon \R \to
[0,1]$ is a smooth monotone increasing function equal to zero near
$-\infty$ and equal to one near $\infty$. Then,
$$
\int_{-\infty}^{\infty}\int_{S^1} \max_W \, 
{\partial}_s H_t^s\,dt\,ds
=\int_{S^1} \max_W \,(H^1_t-H^0_t) \,dt.
$$
Let, as customary, $e^{\scriptscriptstyle{+}}=
e^{\scriptscriptstyle{+}}(H^1-H^0)=\int_{S^1}\max_W \,(H^1_t-H^0_t)
\,dt$. (Recall that $\| H \|= e^{\scriptscriptstyle{+}}(H) -
e^{\scriptscriptstyle{-}}(H)$, where
$e^{\scriptscriptstyle{-}}(H)=\int_{S^1}\min_W \,H_t \,dt$.) Hence, by
Theorem \ref{thm:homotopy}, 
for $a \not\in \Ss(H^s)$, we have the
monotonicity maps $\Psi_{H^0H^1}$ for two intervals: 
$
\HF_*^{(a,\,\infty)} (H^0) \to 
\HF_*^{(a + e^{\scriptscriptstyle{+}},\,\infty )} (H^1)
$
and
$
\HF_*^{(0,\,\infty)} (H^0) \to 
\HF_*^{(e^{\scriptscriptstyle{+}},\,\infty )} (H^1).
$
Now the diagram
\begin{displaymath}
\xymatrix
{
\HF_n^{(0,\,\infty)}(F) \ar[d]_{{\Psi}_{FH^0}} \ar[r]^{{\Psi}_{FH^1}} 
&
\HF_n^{(0,\,\infty)}(H^1) \ar[d]
\\
\HF_n^{(0,\,\infty)}(H^0) \ar[d]_{j_a^{H^0}} \ar[r]^{\Psi_{H^0H^1}}
&
\HF_n^{(e^{\scriptscriptstyle{+}},\,\infty)}(H^1) 
\ar[d]
\\
\HF_n^{(a,\,\infty)}(H^0) \ar[r]^{\Psi_{H^0H^1}}
&
\HF_n^{(a+e^{\scriptscriptstyle{+}},\,\infty)}(H^1),
}
\end{displaymath}
where $F$ is a wide function satisfying $F \ge H^0$ and $F \ge H^1$,
is commutative; see \cite{gi:coiso}. Here the vertical maps on the
right-hand side of the diagram are just the maps induced by taking
quotient complexes, and their composition is the map
$j_{a+e^{\scriptscriptstyle{+}}}^{H^1}$. Consequently, we have
$$
\sigma(H^1) \le \sigma(H^0) + e^{\scriptscriptstyle{+}}
= \sigma(H^0) + e^{\scriptscriptstyle{+}}(H^1-H^0).
$$
Moreover, exchanging the roles of $H^0$ and $H^1$ we get 
$$
\sigma(H^0) \le \sigma(H^1) + e^{\scriptscriptstyle{+}}(H^0-H^1) = 
\sigma(H^1) - e^{\scriptscriptstyle{-}}(H^1-H^0).
$$
Thus, we have 
$$
e^{\scriptscriptstyle{-}}(H^1-H^0) \le \sigma(H^1) - \sigma(H^0) 
\le e^{\scriptscriptstyle{+}}(H^1-H^0).
$$
$C^0$-continuity of the selector follows immediately.  In order to
prove the continuity in Hofer's norm, note first that
$e^{\scriptscriptstyle{+}}(H) \ge 0$ and $e^{\scriptscriptstyle{-}}(H)
\le 0$ for any compactly supported Hamiltonian $H$. (This is also true
for any Hamiltonian on a closed manifold.) Consequently, we have
\begin{eqnarray*}
-e^{\scriptscriptstyle{+}}(H^1-H^0)+e^{\scriptscriptstyle{-}}(H^1-H^0) 
& \le & e^{\scriptscriptstyle{-}}(H^1-H^0) \\
& \le & \sigma(H^1) - \sigma(H^0) \\
& \le & e^{\scriptscriptstyle{+}}(H^1-H^0) \\
& \le & e^{\scriptscriptstyle{+}}(H^1-H^0) - 
e^{\scriptscriptstyle{-}}(H^1-H^0),
\end{eqnarray*}
and hence
$$
|\sigma(H^1) - \sigma(H^0)| 
\le e^{\scriptscriptstyle{+}}(H^1-H^0) - e^{\scriptscriptstyle{-}}(H^1-H^0) 
= \| H^1 -  H^0\|.
$$
This finishes the proof of continuity.

\emph{(AS4)} The assertion readily follows from (AS3).

\emph{(AS5)} We refer the reader to \cite[Lemma 3.5]{gi:weinstein} for
a proof of the property that $\sigma(H) = \max H $ for an autonomous
Hamiltonian $H$ without non-trivial contractible fast periodic orbits.
The proof in \cite{gi:weinstein} is set theoretic in nature and works
in any setting where the selector has the properties (AS0), (AS1),
(AS3) and the claimed property holds for autonomous $C^2$-small
functions. (See the proof of (AS2) above for the fact that $\sigma(H)
= \max H $ when $H$ is a $C^2$-small autonomous function having no
non-trivial contractible fast periodic orbits.)

\emph{(AS6)} It is well-known that the action spectrum of a compactly
supported Hamiltonian on an open symplectically symplectically
aspherical manifold depends only on the time-one flow; see, e.g.,
\cite{FS,hz:book}. Thus, if $H$ and $K$ generate the same time-one
flow and ${\varphi}_H^t$ and ${\varphi}_K^t$ are homotopic (with fixed
end points) in $\Ham_c^+(W)$, then the action spectrum stays the same
throughout the homotopy. On the other hand, due to (AS3), $\sigma$
varies continuously in the course of the homotopy. As the action
spectrum is nowhere dense, $\sigma$ must be constant.

\emph{(AS7)} It is a standard fact that an action selector defined on
a displaceable domain in a closed or convex manifold is \emph{a
priori} bounded from above. However, the proofs existing in literature
rely on the sub-additivity of the action selector and the fact that
the selectors are defined for all Hamiltonians (in particular, not
necessarily non-negative Hamiltonians). Hence, these arguments do not
apply to the action selector introduced here for wide manifolds, and,
for the sake of completeness, we provide a proof of (AS7).

Assume that $V \subset W$ is open and displaceable, and denote by
$e_V$ the displacement energy of $V$; see, e.g., \cite{hz:book, felix,
pol2}. Let $H \colon S^1 \times V \to [0,\,\infty)$ be a compactly
supported Hamiltonian whose support is contained in $V$. Let $K$ be a
compactly supported Hamiltonian such that $\varphi_K$ displaces
$V$. Moreover, without loss of generality, we may assume that $K \ge
0$. For, otherwise, we first shift $K$ up so that $\min K=0$ and then
cut it off away from its original support. The new function is
non-negative, still displaces $V$ and has the same Hofer norm as the
original function.

Recall that, in general, for any two Hamiltonians $H$ and $K$
generating the time-one flows $\varphi_H$ and $\varphi_K$, the
Hamiltonian generating the composition flow $\varphi_H \, \varphi_K$
is given by $H\#K=H(t,x)+K(t,{(\varphi_H^t)}^{-1}(x))$.  Since $K \ge
0$ and the selector is monotone, we then have
\begin{equation}
\labell{eq:sel(H)=sel(HK)}
\sigma(H) \le \sigma(H\#K).
\end{equation}
Since $\varphi_K$ displaces $\supp H$, one-periodic orbits of
$\varphi_H^t \, \varphi_K^t$ are exactly the one-periodic orbits of
$\varphi_K^t$. In fact, we claim that $\Ss(H\#K)=\Ss(K)$. Observe that
this assertion immediately follows from
$$
 \Ss(\varphi_H^t \, \varphi_K^t)
=\Ss(\varphi_K^t  \, \ast \, (\varphi_H^t \, \varphi_K)  )
=\Ss(\varphi_K^t),
$$
where $\ast$ denotes the concatenation of $\varphi_K^t$ and
$\varphi_H^t \, \varphi_K $. Here the concatenation $a(t) \ast b(t)$
of paths $a(t)$ and $b(t)$ with domain $[0,\,1]$ is defined by
traversing $a(2t)$ for $0 \le t \le 1/2$ and then traversing $b(2t-1)$
for $1/2 \le t \le 1$.

The first identity above is due to the fact that $\varphi_H^t \,
\varphi_K^t$ and $\varphi_K^t \, \ast \, (\varphi_H^t \, \varphi_K) $
are homotopic with fixed end points. (It is straightforward to write a
specific formula for this homotopy.)  The second identity is specific
to our situation. Namely, observe that one-periodic orbits of the
concatenation cannot be in $\supp H$, essentially since $\varphi_K$
displaces this support. But when a point is outside $\supp H$, the
flow $\varphi_H^t$ is \emph{identity} and, hence, such a point can
correspond to a one-periodic orbit of the concatenation only if it is
fixed by $\varphi_K$. Therefore, one-periodic orbits of the
concatenation are just reparametrizations of one-periodic orbits of
$\varphi_K^t$.  As a result, actions acquired in both cases are the
same. This proves the second identity.

We now have $\Ss(H\#K)=\Ss(K)$.  The same, of course, holds when $H$
is replaced by $\lambda H$ where $\lambda \in [0,\,\infty)$, i.e.,
$\Ss(\lambda H\#K)=\Ss(K)$ for any non-negative $\lambda$. But
$\sigma$ is continuous and the action spectrum is nowhere
dense. Therefore, we conclude from $\sigma(\lambda H\#K) \in
\Ss(\lambda H\#K)=\Ss(K)$ that $\sigma(\lambda H\#K)$ is independent
of $\lambda$. Setting $\lambda =0$ and $\lambda =1$ yields
$\sigma(K)=\sigma(H\#K)$. Thus, also using \eqref{eq:sel(H)=sel(HK)}
and property (AS4), we have $$ \sigma(H) \le \sigma(H\#K) = \sigma(K)
\le \| K \|.$$

Finally, $e_V = \sup_{K} \| K \|$. Hence, the selector is \emph{a
priori} bounded from above by $e_V$, which is finite when $V$ is
displaceable.


\section{Proofs} 
\labell{proof}

In this section we prove the main results of this paper. 


\subsection{The Conley Conjecture} 
\labell{pf:cc} 
We will focus on the case $W$ is open. For closed manifolds, Theorem
\ref{thm:cc} follows from Theorem \ref{thm:displacement} and the
results from \cite{sc}.

In what follows we denote by $\Ss^+(\cdot)$ the positive part of the
action spectrum of a Hamiltonian.  Theorem \ref{thm:cc} is a
consequence of the following proposition.

\begin{Proposition}
\labell{prop:main}

Let $V$ be an open displaceable subset of a symplectically aspherical
manifold $(W,\omega)$ which is geometrically bounded and wide. Let $G
\geq 0$ be a non-zero Hamiltonian supported in $V$.  Then, $\varphi_G$
has infinitely many periodic points with positive action,
corresponding to contractible periodic orbits of $G$. Moreover, assume
that $\Ss^+(G)$ is separated from zero, i.e., $\mathrm{inf}\,\Ss^+(G)
> 0$. Then, there exists a sequence of integer periods $T_k \to
\infty$ such that for every $T_k$, the Hamiltonian $G$ has such a
periodic orbit with minimal period $T_k$.

\end{Proposition}

Let us first derive the proof of Theorem~\ref{thm:cc} from
Proposition~\ref{prop:main}.

\begin{proof}[Proof of Theorem \ref{thm:cc}.]
Consider $M \times S^1 \subset W \times T^*S^1$, where $T^*S^1$ is
equipped with the standard symplectic structure, also referred to as
the ``stabilization'' of $M$, \cite{mac1, pol1}.  Note that $M \times
S^1$ is again nowhere coisotropic and, moreover, the normal bundle to
$M \times S^1$ admits a non-vanishing section. Theorem
\ref{thm:displacement} now implies that a small neighborhood $V=U
\times (S^1 \times (-\epsilon, \epsilon)) \subset W \times T^*S^1$ of
the product $M \times S^1$ is infinitesimally displaceable. Here $U
\subset W$, a neighborhood of $M$ in $W$, and $\epsilon>0$ are both
sufficiently small.

Let $H \ge 0$ be a Hamiltonian as in the statement of Theorem
\ref{thm:cc}. Let $K \colon T^*S^1 \to [0,1]$ be an autonomous
fiber-wise bump function, depending only on the distance to the zero
section, which is supported in $S^1 \times (-\epsilon, \epsilon)$ and
such that $\max K=K|_{S^1\times 0}=1$. Note that $K$ has no
non-trivial contractible periodic orbits.

Consider the Hamiltonian $G=H \cdot K$ supported in the displaceable
open set $V$. Then $G \geq 0$ and $G \not\equiv 0$. Observe that every
contractible periodic orbit of the vector field $X_{G_t}={H_t} \cdot
X_K+K \cdot X_{H_t}$ with positive action must be of the form
$(u(t),\,v(t)) \in W \times T^*S^1 $, where $u(t)$ is contractible in
$W$ and $v(t)$ is constant with $K(v(t))=1$, i.e., $v(t)$ is a point
on $S^1$. To see this, note first that $v'(t)={H_t}(u(t)) \cdot
X_K(v(t))$ where $X_K$ points in the direction of the angular
coordinate. The assumption that $H\geq 0$ then implies that $v(t) \in
T^*S^1$ can be contractible only when $v(t)$ is constant. Since the
pair $(u(t),\,v(t))$ is contractible, $u(t)$ must also be
contractible. Furthermore, the action on the orbit $(u(t),\,v(t))$ can
be positive only when $v(t)$ is a point on the zero-section.

Coming back to the proof of Theorem \ref{thm:cc}, note that by the
previous observation we have $\Ss^+(H)=\Ss^+(G)$. Moreover, since
$\varphi_H$ is assumed to have isolated fixed points with positive
action, $\Ss^+(H)$ is separated from zero. Then, $\Ss^+(G)$ is also
separated from zero and Proposition \ref{prop:main} applies. Finally,
note that $T_k$-periodic orbits of $\varphi_G$ from the proposition
will correspond to infinitely many simple (contractible) periodic
orbits of $\varphi_H$.
\end{proof}

Let us now prove Proposition \ref{prop:main}.

\begin{proof}[Proof of Proposition \ref{prop:main}.]
Note first that, by property (AS7) and the assumption that $V$ is
displaceable, we have $\sigma(G) \leq e_V < \infty$, where $e_V$ is
the displacement energy of $V$.

Let $G^k$ denote the Hamiltonian generating $\varphi_G^k$. Then, since
$G\geq0$, we have $G^l \leq G^k$ whenever $l<k$. (Here we are using
the explicit formula for the Hamiltonian generating the composition
flow; see, e.g., the proof of (AS7) for this formula.)  By the
monotonicity and non-degeneracy of and the \emph{a priori} bound on
$\sigma$, we then have the following series of inequalities:
$$
0<\sigma(G)\leq\sigma(G^2)\leq \ldots \leq e_V. 
$$

Assume the contrary: $\varphi_G$ has finitely many (simple) periodic
points with positive action. Then, for a sufficiently large prime
number $p$, fixed points of $\varphi_G^p$ must all be $p$-th
iterations of fixed points of $\varphi_G$. Consequently, $\Ss
^+(\varphi_G^p)=p\, \Ss ^+(\varphi_G)$ and, in particular,
$\sigma(G^p) \in p\,\Ss ^+(\varphi_G)$.

Now observe that it suffices to prove the ``moreover'' assertion of
the proposition.  For when $\Ss^+(G)$ is not separated from zero
$\varphi_G$ has infinitely many fixed points and hence has infinitely
many periodic points. Thus, assume that $\inf \Ss^+(G) > \delta > 0 $
for a sufficiently small $\delta > 0$.

Since $\sigma(G^p) >0$ and $\sigma(G^p) \in p\,\Ss ^+(G)$, we
necessarily have $\sigma(G^p) > p \cdot \delta >0$. Hence,
$\sigma(G^p) \to \infty$ as $p \to \infty$. This contradicts the fact
that the selector is \emph{a priori} bounded from above.
\end{proof}

\begin{Remark}
In fact, we have proved that for any prime number $p > e_V / \delta $
there exists a simple contractible $p$-periodic orbit of $G$.
Furthermore, the number of simple (contractible) periodic orbits of
$\varphi_G$ with period less than or equal to $k \in \N$ is at least a
constant, depending on $\varphi_G$, times $k$. This can be seen by
applying the argument in \cite[Proposition 4.13]{vi2}.
\end{Remark}

\subsection{Almost Existence} 
\labell{pf:ae}

We will next prove the almost existence theorem.  

\begin{proof}[Proof of Theorem \ref{thm:ae}.]
The assertion follows from Theorem \ref{thm:displacement} and the
results established in \cite{felix}.  Namely, similarly to the proof
of Theorem \ref{thm:cc}, a sufficiently small neighborhood of $M
\times S^1 \subset P \times T^*S^1$ is infinitesimally displaceable by
the displacement principle for nowhere coisotropic submanifolds. Let
$V=U \times (S^1 \times (-\epsilon, \epsilon)) \subset P \times
T^*S^1$ be such a neighborhood and let, as before, $e_V$ denote the
displacement energy of $V$.

Let us recall the definition of the contractible Hofer-Zehnder
capacity $\Cc(U)$ of a domain $U \subset P $. Denote by
${\Hh}_{\HZ}(U) $ the space of compactly supported smooth non-negative
Hamiltonians $H \colon U \to \R$ which are constant near their maxima
and have no non-trivial \emph{contractible-in-$P$} fast periodic
orbits. The contractible Hofer-Zehnder capacity is then defined to be
$$
\Cc(U)=\sup_H \, \{ \max H \,|\, H \in {\Hh}_{\HZ}(U) \}. 
$$
Removing the condition that the orbits are contractible in $P$ yields
a finer capacity, the ordinary Hofer-Zehnder capacity
$\CHZ(U)$. Hence, $\CHZ(U) \leq \Cc(U)$.

As an immediate consequence of the finiteness of $e_V$ and the
energy-capacity inequality established in \cite[Theorem 1.1]{felix},
we obtain the estimate
$$ 
\CHZ(U) \leq \Cc(U) \leq 4 \, e_V < \infty.
$$ 
In particular, $\CHZ(U) < \infty $ and the almost existence theorem
follows from the finiteness of $\CHZ(U)$ by the standard arguments;
see, e.g., \cite{FGS, hz:book, gg3}.
\end{proof}

\begin{Remark}
As we have already mentioned, the displacement energy-capacity
inequality of \cite{felix} relies heavily on the work of
Lalonde--McDuff and McDuff--Slimowitz,
\cite{LalMc1,McDSl}. Let us give a simple proof of Theorem \ref{thm:ae}
for symplectically aspherical wide manifolds $(W,\omega)$ using the
selector we have constructed in Section \ref{sel-defn}. 

Recall that using an action selector one can define an invariant, the
homological capacity $\chom(V)$, of a domain $V \subset W \times
T^*S^1$ as follows: 
$$
\chom(V)=\sup \{\sigma(H)\,|\, H \colon S^1 \times V \to \R,
\text{ where } H\geq 0 \text{ is compactly supported} \}.
$$
By definition $\chom(V)< \infty$, provided that the selector is
\emph{a priori} bounded from above; for instance, when $V$ is
displaceable. In our case, this is guaranteed by property
(AS7). Furthermore, property (AS5) implies that $\Cc(V) \leq
\chom(V)$. Hence, we obtain $\Cc(V) < \infty$ whenever $V$ is
displaceable, and, consequently, the almost existence theorem holds
for $V$. An argument similar to the one in the proof of Theorem
\ref{thm:cc} finishes the proof.

It is clear that this proof also works for closed manifolds, albeit
using the selector constructed in \cite{sc}.
\end{Remark}

\subsection{Displacement Principle} 
\labell{pf:displacement} 

The proof of Theorem \ref{thm:displacement} is essentially identical
to, if not slightly simpler than, the proof of this statement in the
middle-dimensional case due to Laudenbach and Sikorav,
\cite{LauSi}. Therefore, we will only outline the proof of this
theorem and refer the reader to \cite{LauSi} for the details of the
argument.

Note that the bundle $TM^{\omega}$ is isomorphic to the normal bundle
to $M$ and hence admits a non-vanishing section. To prove Theorem
\ref{thm:displacement}, it suffices to find a non-vanishing section
$v$ of $TM^{\omega}$ such that $dK(v)>0$ along $M$ for some function
$K$ defined near $M$. (The Hamiltonian vector field $X_K$ of $K$ would
then be nowhere tangent to $M$.) On the other hand, for a fixed $v$,
the existence of such a function is guaranteed by a result due to
Sullivan, \cite{su}, whenever $v$ is non-recurrent, i.e., no
trajectory of $v$ is contained entirely in $M$. Then the proof of
Theorem \ref{thm:displacement}, similarly to the argument in
\cite{LauSi}, is based on constructing a non-recurrent section of
$TM^{\omega}$.

Let $\xi$ be a non-vanishing section of $TM^{\omega}$. \emph{A priori}
this section may be somewhere tangent to $M$. The idea is to turn
$\xi$ into a non-recurrent vector field by killing the recurrence. To
this end, let $R \subset M$ denote the set of all trajectories of
$\xi$ contained in $M$. Pick finitely many (mutually) disjoint
balls $B_i \subset M$ in such a way that every trajectory of
$\xi$ in $M$ intersects the interior of at least one $B_i$. This is
possible since $M$ is compact. Denote by $B$ the (disjoint) union of
$B_i$'s.

Next let us modify $\xi$ inside $B$ to make $R=\emptyset$. Observe
that, if the balls are chosen small enough, inside each $B_i$ the
vector field $\xi$ is almost tangent to $M$ by continuity. Thus, let
us choose a \emph{normal} vector $\zeta_i \in TM^{\omega}$ within each
$B_i$ and add these to $\xi$ using cut-off functions supported in
$B_i$'s. The resulting vector field $v$ has the desired property.
(The real situation is slightly more complicated, for actually $\xi$
need not be tangent to $M$ everywhere in $B_i$ and adding $\zeta_i$ to
$\xi$ may force $v$ to vanish in some $B_i$.  However, Thom's jet
transversality theorem guarantees that $\zeta_i$'s can be chosen so
that they not parallel to $\xi$ in each $B_i$.) This finishes the
construction.  Let us point out that the essence of assuming $M$ to be
nowhere coisotropic is now more transparent: in order for $\zeta_i's$
to exist, some ``normal space'' in $TM^{\omega}$ is needed; for
instance, this would be impossible if $TM^{\omega} \subset TM$ at a
point in $R$.

Now it is easy to see that $v$ is non-recurrent, just as in
\cite{LauSi}.



\begin{thebibliography}{CFHWZ}

\bibitem[AL]{al}
M. Audin, J. Lafontaine (Eds), 
\emph{Holomorphic curves in symplectic geometry}, Progress in
Mathematics, \textbf{117}, Birkh\"auser Verlag, Basel, 1994.


\bibitem[BPS]{bps}
P. Biran, L. Polterovich, D. Salamon,
Propagation in Hamiltonian dynamics and relative symplectic homology,
\emph{Duke Math.~J.}, \textbf{119} (2003), 65--118.


\bibitem[CGK]{cgk}
K. Cieliebak, V. Ginzburg, E. Kerman,
Symplectic homology and periodic orbits near symplectic submanifolds, 
\emph{Comment. Math. Helv.}, \textbf{79} (2004), 554--581.


\bibitem[Co]{conley}
C.C. Conley, 
Lecture at the University of Wisconsin, April 6, 1984 


\bibitem[FH]{fh}
A. Floer, H. Hofer,  
Symplectic homology, I. Open sets in $\C^n$, \emph{Math.\ Z.},
\textbf{215} (1994), 37--88.


\bibitem[FHS]{FHS}
A. Floer, H. Hofer, D. Salamon,
Transversality in elliptic Morse theory for the symplectic action,
\emph{Duke Math.\ J.}, \textbf{80} (1995), 251--292.


\bibitem[FrHa]{FrHa} J. Franks, M. Handel, 
Periodic points of Hamiltonian surface diffeomorphisms,
\emph{Geom. Topol.}, \textbf{7} (2003), 713--756.


\bibitem[FGS]{FGS}
U. Frauenfelder, V.L. Ginzburg, F. Schlenk,
Energy capacity inequalities via an action selector, in
\emph{Proceedings on Geometry, Groups, Dynamics and Spectral 
Theory in memory of Robert Brooks.}


\bibitem[FS]{FS}
U. Frauenfelder, F. Schlenk, 
Hamiltonian dynamics on convex symplectic manifolds, \emph{Israel
J.~Math}, \textbf{15} (2006).


\bibitem[Gi1]{gi:survey}
V.L. Ginzburg, 
On closed trajectories of a charge in a magnetic field. An application of 
symplectic geometry, in \emph{Contact and symplectic geometry (Cambridge, 
1994)}, 131--148, Publ.\ Newton Inst., 8, 
Cambridge Univ.\ Press, Cambridge, 1996. 


\bibitem[Gi2]{gi:private}
V.L. Ginzburg, private communication, 2005.


\bibitem[Gi3]{gi:weinstein}
V.L. Ginzburg,
The Weinstein conjecture and theorems of nearby and almost existence in
\emph{The Breadth of Symplectic and Poisson Geometry, 139--172}.
Progr. Math., 232, Birkh\"auser Boston, Boston, MA, 2005.


\bibitem[Gi4]{gi:coiso} 
V.L. Ginzburg, 
Coisotropic Intersections, Preprint 2006, math.SG/0605186; to appear in 
\emph{Duke Math. J.}


\bibitem[Gi5]{gi:conley} 
V.L. Ginzburg, 
The Conley Conjecture, Preprint 2006, math.SG/0610956


\bibitem[GG]{gg3} 
V.L. Ginzburg, B.Z. G\"urel, 
Relative Hofer-Zehnder capacity and periodic orbits in twisted
cotangent bundles, \emph{Duke Math. J.}, \textbf{123} (2004), 1--47.


\bibitem[GK1]{gk1}
V.L. Ginzburg, E. Kerman,
Periodic orbits in magnetic fields in dimensions greater than two,  
in \emph{Geometry and topology in dynamics (Winston-Salem, NC, 
1998/San Antonio, TX, 1999)}, 113--121, Contemp.\ Math., \textbf{246}, 
Amer.\ Math.\ Soc., Providence, RI, 1999. 


\bibitem[GK2]{gk2}
V.L. Ginzburg, E. Kerman,
Periodic orbits of Hamiltonian flows near symplectic extrema,
\emph{Pacific J.\ Math.}, \textbf{206} (2002), 69--91. 


\bibitem[Hi]{hingston}
N. Hingston, 
Subharmonic solutions of Hamiltonian equations on tori, Preprint 2004.


\bibitem[HV]{hv}
H. Hofer, C. Viterbo, 
The Weinstein conjecture in the presence of holomorphic spheres,
\emph{Comm.\ Pure Appl.\ Math.}, \textbf{45} (1992), 583--622.


\bibitem[HZ1]{hz:ae}
H. Hofer, E. Zehnder,
Periodic solutions on hypersurfaces and result by C. Viterbo,
\emph{Invent. Math.} \textbf{90} (1987), 1--9.


\bibitem[HZ2]{ho-ze:capacity}
H. Hofer, E. Zehnder, 
A new capacity for symplectic manifolds, in \emph{Analysis, et
cetera}, P. Rabinowitz and E. Zehnder (Eds.), Academic Press, Boston,
MA, 1990, pp. 405--427.


\bibitem[HZ3]{hz:book}
H. Hofer, E. Zehnder, 
\emph{Symplectic invariants and Hamiltonian dynamics}, Birkh\"{a}user,
Advanced Texts; Basel--Boston--Berlin, 1994.


\bibitem[Ke1]{ke1} 
E. Kerman, 
Periodic orbits of Hamiltonian flows near symplectic critical
submanifolds, \emph{IMRN} (1999), no.\ 17, 954--969.


\bibitem[Ke2]{ke3}
E. Kerman, 
Squeezing in Floer theory and refined Hofer-Zehnder capacities of sets
near symplectic submanifolds, \emph{Geom. Topol.}, \textbf{9} (2005)
1775--1834.


\bibitem[LalMc1]{LalMc1}
F. Lalonde, D. McDuff, 
Hofer's $L^\infty$-geometry: energy and stability of Hamiltonian
flows, part I, \emph{Invent.\ Math.}, \textbf{122} (1995), 1--33.


\bibitem[LalMc2]{LalMc2}
F. Lalonde, D. McDuff, 
Hofer's $L^\infty$-geometry: energy and stability of Hamiltonian
flows, part II, \emph{Invent.\ Math.}, \textbf{122} (1995), 35--69.


\bibitem[LauSi]{LauSi}
F. Laudenbach, J.-C. Sikorav, 
Hamiltonian disjunction and limits of Lagrangian submanifolds,
\emph{IMRN} (1994), no.\ 4, 161--168.


\bibitem[LT]{LT2}
G. Liu, G. Tian, 
Weinstein conjecture and GW invariants, \emph{Commun.\ Contemp.\
Math.}, \textbf{2} (2000), 405--459.

\bibitem[Lu1]{lu}
G. Lu, 
The Weinstein conjecture on some symplectic manifolds containing the
holomorphic spheres, \emph{Kyushu J.\ Math.}, \textbf{52} (1998),
no.\ 2, 331--351. 


\bibitem[Lu2]{lu2}
G. Lu,
Finiteness of Hofer-Zehnder symplectic capacity of neighborhoods of
symplectic submanifolds, Preprint 2005, math.SG/0510172.


\bibitem[Mac]{mac1}
L. Macarini, 
Hofer--Zehnder capacity and Hamiltonian circle actions,
\emph{Commun. Contemp. Math.} \textbf{6} (2004), no. 6, 913--945.


\bibitem[McDSl]{McDSl}
D. McDuff, J. Slimowitz, 
Hofer-Zehnder capacity and length-minimizing Hamiltonian paths,
\emph{Geom. Topol.}, \textbf{5} (2001) 799--830.


\bibitem[Mo]{Mo:orbits}
J. Moser, 
Periodic orbits near equilibrium and a theorem by Alan Weinstein,
\emph{Comm.\ Pure Appl.\ Math.}, \textbf{29} (1976), 727--747.


\bibitem[Oh]{oh3} 
Y.-G. Oh, 
Construction of spectral invariants of Hamiltonian paths on closed
symplectic manifolds, in \emph{The Breadth of Symplectic and Poisson
Geometry, 525--570}. Progr. Math., 232, Birkh\"auser Boston, Boston,
MA, 2005.


\bibitem[Pol1]{pol1}
L. Polterovich, 
An obstacle to non-Lagrangian intersections, in \emph{The Floer
memorial volume}, 575--586, Progr.\ Math., \textbf{133}, Birkh\"auser,
Basel, 1995.


\bibitem[Pol2]{pol2}
L. Polterovich, 
Geometry on the group of Hamiltonian diffeomorphisms.  in
\emph{Proceedings of the International Congress of Mathematicians,
Vol. II (Berlin, 1998).} Doc.\ Math.\ 1998, Extra Vol. II, 401--410
(electronic).


\bibitem[Sa]{sa}
D.A. Salamon, 
Lectures on Floer homology, in \emph{Symplectic Geometry and
Topology}, Eds: Y. Eliashberg and L. Traynor, IAS/Park City
Mathematics series, \textbf{7}, 1999, pp. 143--230.


\bibitem[SZ]{SZ}
D. Salamon, E. Zehnder, 
Morse theory for periodic solutions of Hamiltonian systems and the
Maslov index, \emph{Comm. Pure Appl. Math.},\textbf{45} (1992),
1303--1360.


\bibitem[Schl]{felix}
F. Schlenk, 
Applications of Hofer's geometry to Hamiltonian dynamics,
\emph{Comment. Math. Helv.}, \textbf{81} (2006), 105--121.


\bibitem[Sc]{sc} 
M. Schwarz, 
On the action spectrum for closed symplectically aspherical manifolds,
\emph{Pacific J. Math.}, \textbf{193} (2000), 419--461.


\bibitem[St]{St}
M. Struwe, 
Existence of periodic solutions of Hamiltonian systems on almost every
energy surfaces, \emph{Bol.\ Soc.\ Bras.\ Mat.}, \textbf{20} (1990),
49--58.


\bibitem[Su]{su}
Dennis Sullivan,
Cycles for the dynamical study of foliated manifolds and complex
manifolds, \emph{Invent. Math.}, \textbf{46} (1976), 225--255.


\bibitem[Vi1]{vi2}
C. Viterbo, Symplectic topology as the geometry of generating
functions, \emph{Math.\ Ann.}, \textbf{1992}, 685--710.


\bibitem[Vi2]{vi:functors}
C. Viterbo, 
Functors and computations in Floer homology with applications, I,
\emph{Geom.\ Funct.\ Anal.}, \textbf{9} (1999), 985--1033.



\end{thebibliography}
\end{document}